\documentclass[10pt]{article}
\font\smallit=cmti10
\font\smalltt=cmtt10

\usepackage{amssymb,latexsym,amsmath,epsfig,amsthm} 

\makeatletter

\renewcommand\section{\@startsection {section}{1}{\z@}
{-30pt \@plus -1ex \@minus -.2ex}
{2.3ex \@plus.2ex}
{\normalfont\normalsize\bfseries\boldmath}}

\renewcommand\subsection{\@startsection{subsection}{2}{\z@}
{-3.25ex\@plus -1ex \@minus -.2ex}
{1.5ex \@plus .2ex}
{\normalfont\normalsize\bfseries\boldmath}}

\renewcommand{\@seccntformat}[1]{\csname the#1\endcsname. }

\makeatother

\newcommand{\Z}{\mathbb Z}

\newcommand{\su}{\subseteq}

\newtheorem{theorem}{Theorem}
\newtheorem{lemma}{Lemma}
\newtheorem{proposition}{Proposition}
\newtheorem{corollary}{Corollary}

\theoremstyle{definition}
\newtheorem{definition}{Definition}
\newtheorem{observation}{Observation}

\newtheorem{remark}{Remark}

\begin{document}

\begin{center}
\uppercase{\bf Extremal sequences for a weighted zero-sum constant}
\vskip 20pt
{\bf Santanu Mondal\footnote{Funded by CSIR, Govt. of India.}}\\
{\smallit Department of Mathematics, RKM Vivekananda Educational and Research Institute, P.O. Belur Math, Dist. Howrah, West Bengal, India}\\
{\tt santanu.mondal.math18@gm.rkmvu.ac.in}\\
\vskip 10pt
{\bf Krishnendu Paul}\\
{\smallit Department of Mathematics, RKM  Vivekananda Educational and Research Institute, P.O. Belur Math, Dist. Howrah, West Bengal, India}\\
{\tt krishnendu.p.math18@gm.rkmvu.ac.in}\\
\vskip 10pt
{\bf Shameek Paul}\\
{\smallit Department of Mathematics, RKM  Vivekananda Educational and Research Institute, P.O. Belur Math, Dist. Howrah, West Bengal, India}\\
{\tt shameek.paul@rkmvu.ac.in}\\
\end{center}
\vskip 20pt

\centerline{\smallit Received: , Revised: , Accepted: , Published: } 

\vskip 30pt

\centerline{\bf Abstract}

\noindent
The constant $C_A(n)$ is defined to be the smallest natural number $k$ such that any sequence of $k$ elements in $\Z_n$ has a subsequence of consecutive terms whose $A$-weighted sum is zero, where the weight set $A\su \Z_n\setminus \{0\}$. If $C_A(n)=k$, then a sequence in $\Z_n$ of length $k-1$ which has no $A$-weighted zero-sum subsequence of consecutive terms is called an $A$-extremal sequence. We characterize these sequences for some particular weight sets.

\pagestyle{myheadings}
\markright{\smalltt INTEGERS: 23 (2022)\hfill}
\thispagestyle{empty}
\baselineskip=12.875pt
\vskip 30pt

\section{Introduction}\label{0}

We begin with a very general setup. However, we will soon specialize to the group $\Z_n$, which we will regard as a $\Z_n$-module. Throughout Section \ref{0}, $R$ will be a commutative ring with unity and $A$ will be a non-empty subset of $R$.  

\begin{definition}
Let $M$ be an $R$-module and let $A \su R$. For any $k\geq 1$, a sequence $S=(x_1,\ldots,x_k)$ in $M$ is called an {\it $A$-weighted zero-sum sequence} if for $1\leq i\leq k$, there exists $a_i\in A$ such that $a_1x_1+\cdots+a_kx_k=0$. 
\end{definition}

When $A=\{1\}$, an $A$-weighted zero-sum sequence is  called a zero-sum sequence.

\begin{definition} 
For a finite $R$-module $M$ and for a subset $A \su R$, the {\it $A$-weighted Davenport
constant of $M$} denoted by $D_A(M)$ is defined
to be the least positive integer $k$ such that any sequence in $M$ of length $k$ has an $A$-weighted zero-sum subsequence. 
\end{definition}

This constant was introduced in \cite{AC} for finite abelian groups. It was earlier introduced in \cite{AR} for finite cyclic groups, following a similar generalization in \cite{ACF}. The following related constant was introduced  in $\cite{SKS}$. 

\begin{definition}
For a finite $R$-module $M$ and for a subset $A \su R$, we define the constant $C_A(M)$ to be the least positive integer $k$ such that any sequence in $M$ of length $k$ has an $A$-weighted zero-sum subsequence of consecutive terms. 
\end{definition}

For any $A\su R$, we have $D_A(M)\leq C_A(M)\leq C_{\{1\}}(M)$. The next result shows that for a finite $R$-module $M$, we have  $C_{\{1\}}(M)\leq |M|$. This argument which uses the pigeonhole principle is well-known (see \cite{AC} and \cite{SKS}, for instance). 

\begin{theorem}
Let $M$ be a finite $R$-module with $|M|=n$ and let $S=(x_1,\ldots,x_n)$ be a sequence in $M$ of length $n$. Then there exists $1\leq i\leq j\leq n$ such that $x_i+x_{i+1}+\cdots+x_j=0$.   
\end{theorem}

\begin{proof}
 Let $S=(x_1,\ldots,x_n)$ be a sequence in $M$. For $1\leq i\leq n$, let $y_i=x_1+x_2+\cdots+x_i$. If there exists $i$ such that $y_i=0$, we are done. If $y_i\neq 0$ for all $i$, then there exists $i<j$ such that $y_i=y_j$ and so $x_{i+1}+x_{i+2}+\cdots+x_j=y_j-y_i=0$. 
\end{proof}

\begin{definition}
Let $A\su R$ and let $k=C_A(M)$. Suppose $k\geq 2$. Then  there exists a sequence $S$ in $M$ of length $k-1$ which has no $A$-weighted zero-sum subsequence of consecutive terms. We will call such a sequence an {\it $A$-extremal sequence}. 
\end{definition}

We consider the ring $\Z_n$ as a module over itself and for $A\su\Z_n$, we let $C_A(n)$ denote $C_A(\Z_n)$. We observe that when $n\geq 2$ and $0\notin A$, we have $C_A(n)\geq 2$. In this paper we study the $A$-extremal sequences for the $\Z_n$-module $\Z_n$. Such sequences were studied in \cite{AMP} in the context of the weighted Davenport constant. 

\begin{remark}
If we permute the terms of an $A$-extremal sequence, then we may not get an $A$-extremal sequence. For example, the sequence $S=(2,1,2)$ in $\Z_4$ does not have any $A$-weighted zero-sum subsequence of consecutive terms where $A=\{1,-1\}$. Also, from Corollary 3 of $\cite{SKS}$ we have $C_A(4)=4$ and so $S$ is an $A$-extremal sequence. Clearly $(1,2,2)$ is not an $A$-extremal sequence in $\Z_4$. 
\end{remark}

Let $A=\{1\}$. In Corollary 1 of $\cite{SKS}$ it was shown that $C_A(n)=n$. The next theorem characterizes the $A$-extremal sequences in $\Z_n$. 

\begin{theorem}
Let $A=\{1\}$ and let $S=(x_1,\,\ldots,\,x_{n-1})$ be a sequence in $\Z_n$. For $1\leq i\leq n-1$, let $y_i=x_1+\cdots+x_i$. Then $S$ is an $A$-extremal sequence if and only if the set $\{\,y_i:1\leq i\leq n-1\,\}$ is the same as $\Z_n\setminus\{0\}$. 
\end{theorem}

\begin{proof}
Let $A=\{1\}$ and let $S=(x_1,\,\ldots,\,x_{n-1})$ be a sequence in $\Z_n$. As $C_A(n)=n$, it follows that $S$ will be a $A$-extremal sequence if and only if $S$ does not have any zero-sum subsequence of consecutive terms. Let $i$ and $j$ be such that $1\leq i<j\leq n-1$. Then we have $y_i=y_j$ if and only if $x_{i+1}+\cdots+x_j=0$. As a result, the set $\{\,y_i:1\leq i\leq n-1\,\}$ has $n-1$ elements.
\end{proof}

\begin{remark}
For $A=\{1\}$, to construct an $A$-extremal sequence in $\Z_n$, we can choose $x_1\in\Z_n\setminus\{0\}$ in $n-1$ ways. Given $x_1$, we can choose $x_2$ in $n-2$ ways so that $y_2\in\Z_n\setminus\{0,\,y_1\}$. For $3\leq i\leq n-1$, given $x_1,x_2,\ldots,x_{i-1}$, we can choose $x_i$ in $n-i$ ways so that $y_i\in \Z_n\setminus\{0,\,y_1,\,\ldots,\,y_{i-1}\}$. So there are $(n-1)!$ sequences in $\Z_n$ which are $A$-extremal sequences.
\end{remark}

In Theorem 2 of \cite{SKS} it was shown that for $A=\Z_n\setminus \{0\}$, we have $C_A(n)=2$.

\begin{theorem}\label{nz}
For $A=\Z_n\setminus \{0\}$, a sequence $(x)$ of length 1 is an $A$-extremal sequence if and only if $x$ is a unit. 
\end{theorem}

\begin{proof}
We observe that a sequence $(x)$ of length 1 is an $A$-weighted zero-sum sequence if and only if $x$ is a zero-divisor. 
\end{proof}

Let $U(n)$ denote the multiplicative group of units in the ring $\Z_n$. For $j\geq 1$, let $U(n)^j$ denote the set $\{\,x^j\mid x\in U(n)\,\}$. If $n=p_1p_2\ldots p_k$ where $p_i$ is a prime for $1\leq i\leq k$, then we define $\Omega(n)=k$.

\begin{remark}
When $p$ is a prime, Theorem \ref{nz} characterizes the $U(p)$-extremal sequences in $\Z_p$. 
\end{remark}

In this article, we characterize the $A$-extremal sequences in $\Z_n$ for the following weight sets: 

\begin{itemize}
 \item $A=U(p)^2$ where $p$ is a prime; 
 \item $A=U(p)^3$ where $p$ is a prime;  
 \item $A=U(n)$ where $n$ is odd;  
 \item $A=U(n)^2$ where every prime divisor of $n$ is at least 7; 
 \item $A=U(n)^3$ where $n$ is squarefree and is not divisible by 2, 7 or 13.
\end{itemize}

\section{The case when $A=U(p)^2$ where $p$ is a prime}\label{q}

When $A=U(2)^j$ where $j\geq 1$, Theorem \ref{nz} characterizes the $A$-extremal sequences in $\Z_2$, since  for any $j\geq 1$, we have $U(2)^j=U(2)$.

\begin{theorem}
Let $F$ be a field and let $A$ be a subgroup of $F^*:=F\setminus\{0\}$. A sequence $S=(x,y)$ in $F^*$ has an $A$-weighted zero-sum subsequence if and only if $S$ is an $A$-weighted zero-sum sequence if and only if $x$ and $-y$ are in the same coset of $A$.  
\end{theorem}

\begin{proof}
Let $S=(x,y)$ be a sequence in $F^*$. Observe that $x$ and $-y$ are in the same coset of $A$ if and only if there exists $c\in A$ such that $x=-cy$. This is if and only if there exists $c\in A$ such that $x+cy=0$, which happens if and only if there exists $a,b\in A$ such that $ax+by=0$.
\end{proof}

\begin{corollary}\label{field}
Let $F$ be a field and let $A$ be a subgroup of $F^*$. A sequence $S=(x,y)$ in $F$ does not have an $A$-weighted zero-sum subsequence if and only if $x,y\in F^*$ and $x,-y$ are in different cosets of $A$ in $F^*$. 
\end{corollary}

For an odd prime $p$ we denote $U(p)^2$ by $Q_p$. In Theorem 4 of $\cite{SKS}$ it was shown that when $A=Q_p$ where $p$ is a prime, we have $C_A(2)=2$ and $C_A(p)=3$ when $p\neq 2$.

\begin{corollary}\label{q1}
Let $A=Q_p$ where $p$ is an odd prime. Let $S=(x,y)$ be a sequence in $\Z_p$. Then $S$ is an $A$-extremal sequence in $\Z_p$ if and only if $x,y\in U(p)$ and $x,-y$ are in different cosets of $A$ in $U(p)$.   
\end{corollary}

\begin{proof}
This follows from Corollary \ref{field} by considering the field $\Z_p$ and the subgroup $Q_p$ of $\Z_p^*$.  
\end{proof}

\section{The case when $A=U(p)^3$ where $p$ is a prime}

Let $p$ be a prime. In Section 3 of \cite{SKS}, it was shown that when $p\not\equiv 1~(\textrm{mod}~3)$ we have $U(p)^3=U(p)$. In Theorem 7 and Lemma 2 of \cite{SKS}, it was shown that when $p\equiv 1~(\textrm{mod}~3)$, we have $C_A(p)=3$ when $p\neq 7$ and  $C_A(7)=4$.

When $p\not\equiv 1~(\textrm{mod}~3)$ we have $A=\Z_p\setminus\{0\}$, and so Theorem \ref{nz} characterizes the $A$-extremal sequences in $\Z_p$. 
The next two results characterize the $A$-extremal sequences in $\Z_p$ when $p\equiv 1~(\textrm{mod}~3)$. 

\begin{corollary}\label{cubpx}
Let $A=U(p)^3$ where $p$ is a prime such that $p\equiv 1~(\textrm{mod}~3)$ and $p\neq 7$. Let $S=(x,y)$ be a sequence in $\Z_p$. Then $S$ is an $A$-extremal sequence in $\Z_p$ if and only if $x,y\in U(p)$ and $x,-y$ are in different cosets of $A$ in $U(p)$.  
\end{corollary}

\begin{proof}
Let $S=(x,y)$ be a sequence in $\Z_p$. As $C_A(p)=3$, it follows that $S$ is an $A$-extremal sequence in $\Z_p$ if and only if $S$ has no $A$-weighted zero-sum subsequence. Now from Corollary \ref{field}, upon considering $\Z_p$ as a field, the result follows.  
\end{proof}

The following definition is similar to a definition which was made in \cite{AMP}.

\begin{definition}\label{eq}
Let $S=(x_1,\ldots,x_k)$ and $T=(y_1,\ldots,y_k)$ be two sequences in $\Z_n$ and let $A$ be a subgroup of $U(n)$. We say that $S$ is {\it $A$-equivalent} to $T$ if there exists $c\in U(n)$ and if for $1\leq i\leq k$ there exists $a_i\in A$ such that $cy_i=a_ix_i$. 
\end{definition}

Observe that if $S$ is an $A$-weighted zero-sum sequence, then so is any sequence which is $A$-equivalent to $S$. Thus, if $S$ is an $A$-extremal sequence, then so is any sequence which is $A$-equivalent to $S$.

\begin{proposition}\label{ob}
Let $A=U(7)^3$. Then $S$ is an $A$-extremal sequence in $\Z_7$ if and only if $S$ is $A$-equivalent to $(1,3,1)$. 
\end{proposition}

\begin{proof}
Let $S$ be a sequence in $\Z_7$. Suppose $S$ is an $A$-extremal sequence. As $C_A(7)=4$, it follows that $S$ has length 3. As $S$ has no $A$-weighted zero-sum subsequence of consecutive terms, all terms of $S$ are non-zero. Thus, the terms of $S$ are in the set $\{\pm 1,\pm 2,\pm 3\}$. As $A=\{1,-1\}$, no two consecutive terms of $S$ can be equal up to sign and $S$ cannot be $A$-equivalent to $(1,2,1)$ or to any permutation of $(1,2,3)$. By multiplying by a unit, we see that $S$ is $A$-equivalent to a sequence whose first term is 1. Hence, $S$ has to be $A$-equivalent to a sequence of the form $(1,3,1)$. 

Conversely, as $A=\{1,-1\}$, it is easy to see that any sequence which is $A$-equivalent to $(1,3,1)$ is an $A$-extremal sequence as this sequence does not have any $A$-weighted zero-sum subsequence of consecutive terms. 
\end{proof}

\begin{remark}
Theorem \ref{nz}, Corollary \ref{cubpx} and Proposition \ref{ob} together characterize the $U(p)^3$-extremal sequences in $\Z_p$ when $p$ is a prime.  
\end{remark}

\section{The case when $A=U(n)$} 

Let $m$ be a divisor of $n$. We will refer to the homomorphism $\varphi:\Z_n\to \Z_m$ given by $a+n\Z\mapsto a+m\Z$ as the natural map.

We use the notation $p^r\mid\mid n$ to mean $p^r\mid n$ and  $p^{r+1}\nmid n$. For a prime divisor $p$ of $n$ where $p^r\mid\mid n$, let $A_p=U(p^{r})$ and let $S^{(p)}$ denote the sequence in $\Z_{p^r}$ which is the image of $S$ under the natural map $\Z_n\to \Z_{p^r}$. The following was pointed out in Observation 2.2 in \cite{sg}. We restate it here using our notation.

\begin{observation}[\cite{sg}]\label{observation}
Let $A=U(n)$ and let $S$ be a sequence in $\Z_n$. Then $S$ is an $A$-weighted zero-sum sequence in $\Z_n$ if and only if for every prime divisor $p$ of $n$, $S^{(p)}$ is an $A_p$-weighted zero-sum sequence in $\Z_{p^r}$ where $p^r\mid\mid n$.
\end{observation}

In fact we have the following more general result whose proof is similar to the proof of the previous result. Using the notations as in the second paragraph of this section, we define $A_p^j=\{\,x^j : x\in A_p\,\}$. 

\begin{observation}\label{observation2}
Let $A=U(n)^j$ where $j\geq 1$ and let $S$ be a sequence in $\Z_n$. Then $S$ is an $A$-weighted zero-sum sequence in $\Z_n$ if and only if for every prime divisor $p$ of $n$, we have $S^{(p)}$ is an $A_p^j$-weighted zero-sum sequence in $\Z_{p^r}$ where $p^r\mid\mid n$.  
\end{observation}

The following result is Lemma 5 in \cite{SKS}.

\begin{lemma}[\cite{SKS}]\label{b}
Let $S$ be a sequence in $\Z_n$ and let $p$ be a prime divisor of $n$ which divides every element of $S$. Let $n'=n/p$, $A=U(n)$ and $A'=U(n')$. Then $S$ is an $A$-weighted zero-sum  sequence if $S'$ is an $A'$-weighted zero-sum sequence, where $S'$ is the sequence in $\Z_{n'}$ whose terms are obtained by dividing the terms of $S$ by $p$. 
\end{lemma}

The following result is Lemma 2.1 (ii) in \cite{sg} which we restate here using our terminology.

\begin{lemma}[\cite{sg}]\label{gri}
Let $A=U(p^r)$ where $p$ is an odd prime and $r\geq 1$. If a sequence $S$ over $\Z_{p^r}$ has at least two terms coprime to $p$, then $S$ is an $A$-weighted zero-sum sequence. 
\end{lemma}

In Corollary 4 of $\cite{SKS}$ it was shown that for $A=U(n)$ where $n$ is odd we have $C_A(n)=2^{\Omega(n)}$. In the next theorem we give a method to construct $A$-extremal sequences in $\Z_n$.  

\begin{theorem}\label{unit1'}
Let $A=U(n)$ where $n$ is odd, let $p$ be a prime divisor of $n$, $n'=n/p$ and $A'=U(n')$. Then $S=(p\,u_1,\,\ldots,\,p\,u_k,\, x^*,\,p\,v_1,\,\ldots,\,p\,v_k)$ is an $A$-extremal sequence in $\Z_n$ where $S_1'=(u_1,\,\ldots,\,u_k)$ and $S_2'=(v_1,\,\ldots,\,v_k)$ are $A'$-extremal sequences in $\Z_{n'}$, and $x^*\in\Z_n$ is an element which is not divisible by $p$.  
\end{theorem}

\begin{proof}
As $n'$ is odd, $C_{A'}(n')=2^{\Omega(n')}$. As $S_1'$ is an $A'$-extremal sequence in $\Z_{n'}$ of length $k$, it follows that $k=2^{\Omega(n')}-1$. Let $S$ be defined as in the statement of the theorem. We claim that $S$ does not have any $A$-weighted zero-sum subsequence of consecutive terms. 
Suppose our claim is false. Let $T$ be an $A$-weighted zero-sum subsequence of consecutive terms of $S$. Define sequences $S_1$ and $S_2$ in $\Z_n$ as $S_1=(p\,u_1,\,\ldots,\,p\,u_k)$ and $S_2=(p\,v_1,\,\ldots,\,p\,v_k)$.

Suppose $x^*$ is not a term of $T$. Hence, $T$ is a subsequence of either $S_1$ or $S_2$. Let $T'$ be the sequence in $\Z_{n'}$ whose terms are obtained by dividing the terms of $T$ by $p$. The image of $A$ under the natural map $\Z_n\to\Z_{n'}$ is contained in $A'$. As a result, by dividing the $A$-weighted zero-sum which is obtained from $T$ by $p$, we will get the contradiction that either $S_1'$ or $S_2'$ has an $A'$-weighted zero-sum subsequence of consecutive terms, namely $T'$.

So $x^*$ is a term of $T$. Now considering the $A$-weighted zero-sum which is obtained from $T$, we see that a multiple of $p$ is a unit multiple of $x^*$. This gives the contradiction that  $x^*$ is divisible by $p$. Thus our claim is true. As $n$ is odd, $C_A(n)=2^{\Omega(n)}$. Also as $S$ is a sequence of length $2k+1=2^{\Omega(n)}-1$ in $\Z_n$, hence $S$ is an $A$-extremal sequence in $\Z_n$. 
\end{proof}

Let us consider an example to see how we can use Theorem \ref{unit1'} to construct $A$-extremal sequences. By Theorem \ref{nz} we have the $U(5)$-extremal sequences $(2)$ and $(4)$ in $\Z_5$. Now by Theorem \ref{unit1'} we get the $U(25)$-extremal sequences $(10,4,20)$ and $(10,21,20)$ in $\Z_{25}$. Repeating this process gives us the $U(75)$-extremal sequence $(30,12,60,38,30,63,60)$ in $\Z_{75}$.

\begin{lemma}\label{unit4'}
Let $A=U(n)$ where $n$ is odd, $\ell=2^{\Omega(n)}-1$ and let $q$ be a prime divisor of $n$. Then a sequence $S$ in $\Z_n$ has an $A$-weighted zero-sum subsequence of consecutive terms if $S$ has a subsequence $T$ of consecutive terms of length at least $(\ell+1)/2$, such that each term of $T$ is divisible by $q$.  
\end{lemma}

\begin{proof}
Let $n'=n/q$ and let $A'=U(n')$. As $\ell=2^{\Omega(n)}-1$, it follows that $(\ell+1)/2=2^{\Omega(n)-1}$ and hence $T$ has length at least $2^{\Omega(n')}$. Let $T'$ denote the sequence in $\Z_{n'}$ whose terms are obtained by dividing the terms of $T$ by $q$. As $n$ is odd, we have $n'$ is odd. Now as $T'$ is a sequence of length at least $2^{\Omega(n')}$ in $\Z_{n'}$ and as $C_{A'}(n')=2^{\Omega(n')}$, it follows that $T'$ has an $A'$-weighted zero-sum subsequence of consecutive terms. Now by Lemma \ref{b} we see that $T$ has an $A$-weighted zero-sum  subsequence of consecutive terms. As $T$ is a subsequence of consecutive terms of $S$, it follows that $S$ has an $A$-weighted zero-sum subsequence of consecutive terms. 
\end{proof}

\begin{corollary}\label{unit4}
Let $A=U(n)$ where $n$ is odd. Suppose $S$ is an $A$-extremal sequence in $\Z_n$ and $q$ is a prime divisor of $n$. Then $q$ is coprime to at least one term of $S$. 
\end{corollary}

\begin{proof}
As $S$ is an $A$-extremal sequence in $\Z_n$, if $S$ has  length $\ell$ then $\ell=2^{\Omega(n)}-1$. Suppose no term of $S$ is coprime to $q$. As $\ell$ is at least $(\ell+1)/2$, by Lemma \ref{unit4'} we see that $S$ has an $A$-weighted zero-sum subsequence of consecutive terms. This  contradicts our assumption. Hence, at least one term of $S$ is coprime to $q$. 
\end{proof}

Theorem \ref{unit3} shows that the method in Theorem \ref{unit1'} is the only way to obtain $A$-extremal sequences in $\Z_n$ when $n$ is odd. 

\begin{theorem}\label{unit3}
Let $A=U(n)$ where $n$ is odd. Suppose $S=(x_1,\ldots,x_\ell)$ is an $A$-extremal sequence in $\Z_n$. Then there is a prime divisor $p$ of $n$ such that $p$ divides all the terms of $S$ except the `middle' term $x_{k+1}$ where $k+1=(\ell+1)/2$. Moreover, if $S_1=(x_1,\ldots,x_k)$, $S_2=(x_{k+2},\ldots,x_\ell)$, $n'=n/p$ and $A'=U(n')$, then $S_1'$ and $S_2'$ are $A'$-extremal sequences in $\Z_{n'}$ where $S_1'$ and $S_2'$ denote the sequences in $\Z_{n'}$ which are obtained by dividing the terms of $S_1$ and $S_2$ by $p$. 
\end{theorem}

\begin{proof}
As $S$ is an $A$-extremal sequence in $\Z_n$, it follows that $\ell=2^{\Omega(n)}-1$. Suppose for every prime divisor $q$ of $n$, there are at least two terms of $S$ which are not divisible by $q$. Let $q$ be a prime divisor of $n$ and let $q^r\mid\mid n$. As $n$ is odd, $q$ is odd. Let $S^{(q)}$ be defined as at the beginning of this section. As the sequence $S^{(q)}$ contains at least two units, by Lemma \ref{gri} we see that $S^{(q)}$ is an $A_q$-weighted zero-sum sequence where $A_q=U(q^r)$. As this is true for every prime divisor $q$ of $n$, by Observation \ref{observation} we see that $S$ is an $A$-weighted zero-sum sequence. This contradicts our assumption about $S$. Thus, there is a prime divisor $p$ of $n$ such that at most one term of $S$ is not divisible by $p$. Now by Corollary \ref{unit4} we see that there is exactly one term of $S$, say $x^*$, which is not divisible by $p$.

Suppose $x^*\neq x_{k+1}$ where $k+1=(\ell+1)/2$. Thus, there is a subsequence $T$ of consecutive terms of $S$ of length at least $k+1$ such that $p$ divides every term of $T$. As $T$ has length at least $(\ell+1)/2$, by Lemma \ref{unit4'} we see that $S$ has an $A$-weighted zero-sum subsequence of consecutive terms. This contradicts our assumption about $S$. As a result, $x^*=x_{k+1}$.

Let $n'=n/p$ and for $i=1,2$, let $S_i$ and $S_i'$ be defined as in the statement of the theorem. If $S_1'$ has an $A'$-weighted zero-sum subsequence of consecutive terms, then by Lemma \ref{b} we see that $S_1$ has an $A$-weighted zero-sum subsequence of consecutive terms. As $S_1$ is a subsequence of consecutive terms of $S$, it follows that $S$ has an $A$-weighted zero-sum subsequence of consecutive terms. This contradicts our assumption about $S$. Now as $k+1=(\ell+1)/2=2^{\Omega(n)-1}$, it follows that $S_1'$ has length $2^{\Omega(n')}-1$ and hence $S_1'$ is an $A'$-extremal sequence in $\Z_{n'}$. A similar argument shows that $S_2'$ is also an $A'$-extremal sequence in $\Z_{n'}$.
\end{proof}

\begin{remark} 
Theorems \ref{unit1'} and \ref{unit3} together characterize the $U(n)$-extremal sequences in $\Z_n$ where $n$ is odd. 
\end{remark}

\begin{proposition}
Let $A=U(n)$ where $n$ is an odd integer and let $n=p_1p_2$ be a product of two primes which are not necessarily distinct. Then $S$ is an $A$-extremal sequence in $\Z_n$ if and only if $S$ is of the form $(b_1q_1,\,a_1,\,b_2q_1)$ where $q_1\nmid a_1$, $q_2\nmid b_1b_2$ and $q_1,q_2$ is a permutation of $p_1,p_2$. 
\end{proposition}

\begin{proof}
Suppose $S$ is an $A$-extremal sequence in $\Z_n$. As $\Omega(n)=2$, it follows that $S$ has length  $2^{\Omega(n)}-1=3$. Let $S$ be the sequence $(x_1,x_2,x_3)$. By Theorem \ref{unit3} there is a prime $q_1\in\{p_1,p_2\}$ such that $x_2$ is the only term of $S$ which is coprime to $q_1$. Let $n'=n/q_1$, $x_1'=x_1/q_1$, $x_3'=x_3/q_1$ and $A'= U(n')$. Again by Theorem \ref{unit3} we get that the sequences $(x_1')$ and $(x_3')$ are $A'$-extremal sequences in $\Z_{n'}$. As $n'$ is a prime, say $q_2$, it follows that $x_1',x_3'\in U(q_2)$. Thus, $S$ is of the form $(b_1q_1,a_1,b_2q_1)$ where $q_1\nmid a_1$ and $q_2\nmid b_1b_2$ and where  $q_1,q_2$ is a permutation of $p_1,p_2$. The converse of this statement follows from Theorem \ref{unit1'}.
\end{proof}

By a similar argument using Theorem \ref{unit3} as in the previous proof, we get the next result. 

\begin{proposition}
Let $A=U(n)$ where $n$ is an odd integer and let $n=p_1p_2p_3$ be a product of three primes which are not necessarily distinct. Then $S$ is an $A$-extremal sequence in $\Z_n$ if and only if $S$ has one of the following two possible forms. In these two sequences $q_1,q_2,q_3$ is a permutation of $p_1,p_2,p_3$ and for all $i$ we have $a_i,b_i,c_i\in\Z_n$. \\ 

Type (1): $(a_1q_1q_2,\,b_1q_1,\,a_2q_1q_2,\,c_1,\,a_3q_1q_2,\,b_2q_1,\,a_4q_1q_2)$ 
where $q_1\nmid c_1$, $q_2\nmid b_1b_2$ and $q_3\nmid a_1a_2a_3a_4$ \\ 
    
Type (2): $(a_1q_1q_2,\,b_1q_1,\,a_2q_1q_2,\,c_1,\,b_2q_1q_3,\,a_3q_1,\,b_3q_1q_3)$ 
where $q_1\nmid c_1$, $q_2\nmid b_1b_2b_3$ and $q_3\nmid a_1a_2a_3$
\end{proposition}

\section{The case when $A=U(n)^2$}

We need the following result which is Lemma 1 in \cite{CM}. 

\begin{lemma}[\cite{CM}]\label{cm'}
Let $A=U(n)^2$ where $n=p^r$ and $p\geq 7$ is a prime. Let $x_1,x_2,x_3\in U(n)$. Then we have  $Ax_1+Ax_2+Ax_3=\Z_n$. 
\end{lemma}

The next result is Lemma 7 of \cite{SKS} which is an easy consequence of Lemma \ref{cm'}. 

\begin{corollary}[\cite{SKS}]\label{cm}
Let $A=U(n)^2$ where $n=p^r$ and $p\geq 7$ is a prime. Then a sequence $S$ in $\Z_n$ is an $A$-weighted zero-sum sequence if at least three terms of $S$ are in $U(n)$. 
\end{corollary}

\begin{remark}
The conclusion of Corollary \ref{cm} may not hold when $p<7$. The sequence $(1,1,1)$ in $\Z_n$ is not a $U(n)^2$-weighted zero-sum sequence when $n=2,5$. The sequence $(1,2,1)$ in $\Z_3$ is not a $U(3)^2$-weighted zero-sum sequence. 
\end{remark}

The next result is in Lemma 5 of \cite{SKS}.

\begin{lemma}[\cite{SKS}]\label{c}
Let $S$ be a sequence in $\Z_n$, let $p$ be a prime divisor of $n$ which divides every element of $S$. Let  $n'=n/p$, $A=U(n)^2$ and $A'=U(n')^2$. Then $S$ is an $A$-weighted zero-sum sequence if $S'$ is an $A'$-weighted zero-sum sequence, where  $S'$ is the sequence in $\Z_{n'}$ whose terms are obtained by dividing the terms of $S$ by $p$.   
\end{lemma}

In Corollary 6 of $\cite{SKS}$ it was shown that for $A=U(n)^2$ where every prime divisor of $n$ is at least 7, we have $C_A(n)=3^{\Omega(n)}$. We now give a method to construct $A$-extremal sequences in $\Z_n$.

\begin{theorem}\label{sq}
Let $A=U(n)^2$ where every prime divisor of $n$ is at least 7, and let \newline $p$ be a prime divisor of $n$. Let $n'=n/p$ and $A'=U(n')^2$. Then the sequence $S=(p\,u_1,\,\ldots,\,p\,u_k,\, x^*,\,p\,v_1,\,\ldots,\,p\,v_k,\,x^{**},\, p\,w_1,\,\ldots,\,p\,w_k)$ is an $A$-extremal sequence in $\Z_n$ where $S_1'=(u_1,\,\ldots,\,u_k)$, $S_2'=(v_1,\,\ldots,\,v_k)$ and $S_3'=(w_1,\,\ldots,\,w_k)$ are $A'$-extremal sequences in $\Z_{n'}$  and $x^*,x^{**}\in\Z_n$ are such that the image of the sequence $(x^*,x^{**})$ under the natural map $\Z_n\to\Z_p$ does not have any $Q_p$-weighted zero-sum subsequence.   
\end{theorem}

\begin{proof}
Since $S_1'$ is an $A'$-extremal sequence of length $k$ in $\Z_{n'}$, so $k=3^{\Omega(n')}-1$.
Let $S$ be defined as in the statement of the theorem. We claim that $S$ does not have any $A$-weighted zero-sum subsequence of consecutive terms. Suppose our claim is false. Let $T$ be an $A$-weighted zero-sum subsequence of consecutive terms of $S$. Define sequences $S_1,S_2$ and $S_3$ in $\Z_n$ as $S_1=(p\,u_1,\,\ldots,\,p\,u_k),~S_2=(p\,v_1,\,\ldots,\,p\,v_k)$, $S_3=(p\,w_1,\,\ldots,\,p\,w_k)$.

The image of $A$ under the natural map $\Z_n\to\Z_{n'}$ is contained in $A'$. So if $T$ is a subsequence of either $S_1$, $S_2$ or $S_3$, then as we have seen in the second paragraph of the proof of Theorem \ref{unit1'}, we will get the contradiction that either $S_1'$, $S_2'$ or $S_3'$ has an $A'$-weighted zero-sum subsequence of consecutive terms. So $T\cap \{x^*,x^{**}\}\neq\varnothing$. Now considering the $A$-weighted zero-sum which is obtained from $T$, we see that an $A$-weighted linear combination of a subsequence of  $(x^*,x^{**})$ is a multiple of $p$.

The image of $A$ under the natural map $\Z_n\to\Z_p$ is contained in $Q_p$. So we get a contradiction that the image of the sequence $(x^*,x^{**})$ under the natural map $\Z_n\to\Z_p$ has a $Q_p$-weighted zero-sum subsequence. Thus our claim is true. Now $S$ is a sequence of length $3k+2=3^{\Omega(n)}-1$ in $\Z_n$. Hence, $S$ is an $A$-extremal sequence in $\Z_n$. 
\end{proof}

We can use Theorem \ref{sq} to construct $A$-extremal sequences in $\Z_n$. For example, from Corollary \ref{q1} we get the $U(7)^2$-extremal sequences $(4,1)$, $(3,3)$, $(2,1)$ and $(1,1)$ in $\Z_7$. Now by using these we get $(28,7,15,21,21,36,14,7)$ which is a $U(49)^2$-extremal sequence in $\Z_{49}$.

Theorems \ref{7} and \ref{8} will show us that the method in Theorem \ref{sq} is the only way to obtain $A$-extremal sequences in $\Z_n$ when any prime divisor of $n$ is at least $7$. We shall need the following lemma.

\begin{lemma}\label{g'}
Let $A=U(n)^2$ where every prime divisor of $n$ is at least 7, $\ell=3^{\Omega(n)}-1$ and let $q$ be a prime divisor of $n$. Then a sequence $S$ in $\Z_n$ has an $A$-weighted zero-sum subsequence of consecutive terms if $S$ has a subsequence $T$ of consecutive terms whose length is at least $(\ell+1)/3$, such that each term of $T$ is divisible by $q$. 
\end{lemma}

\begin{proof}
Let $n'=n/q$ and let $T'$ denote the sequence in $\Z_{n'}$ whose terms are obtained by dividing the terms of $T$ by $q$. As $T$ is a sequence of length at least $(\ell+1)/3$ in $\Z_n$, it follows that $T'$ is a sequence of length at least $3^{\Omega(n')}$ in $\Z_{n'}$. Also, as every prime divisor of $n'$ is at least 7 and as $C_{A'}(n')=3^{\Omega(n')}$, it follows that $T'$ has an $A'$-weighted zero-sum subsequence of consecutive terms where $A'=U(n')^2$. Now by Lemma \ref{c} we see that $T$ has an $A$-weighted zero-sum subsequence of consecutive terms. As $T$ is a subsequence of consecutive terms of $S$, it follows that $S$ has an $A$-weighted zero-sum subsequence of consecutive terms. 
\end{proof}

\begin{corollary}\label{g}
Let $A=U(n)^2$ where every prime divisor of $n$ is at least 7. Suppose $S$ is an $A$-extremal sequence in $\Z_n$. Then for any prime divisor $q$ of $n$, we have that $q$ is coprime to at least two terms of $S$. 
\end{corollary}

\begin{proof}
As $S$ is an $A$-extremal sequence in $\Z_n$, so $\ell=3^{\Omega(n)}-1$ if $S$ has length $\ell$. Suppose there is a prime divisor $q$ of $n$ such that at most one term of $S$ is coprime to $q$. Then we get a subsequence, say $T$, of consecutive terms of $S$ of length at least $(\ell+1)/3$ such that each term of $T$ is divisible by $q$. Now by Lemma \ref{g'} we see that $S$ has an $A$-weighted zero-sum subsequence of consecutive terms. This contradicts our assumption about $S$. Hence, at least two terms of $S$ are coprime to $q$. 
\end{proof}

\begin{theorem}\label{7}
Let $A=U(n)^2$ where any prime divisor of $n$ is at least $7$. Suppose $S=(x_1,\ldots,x_\ell)$ is an $A$-extremal sequence in $\Z_n$. Then there is a prime divisor $p$ of $n$ such that $p$ divides all the terms of $S$ except the terms $x_{k+1}$ and $x_{2(k+1)}$ where $k+1=(\ell+1)/3$.
Moreover, if $S_1=(x_1,\ldots,x_k)$, $S_2=(x_{k+2},\ldots,x_{2k+1})$, $S_3=(x_{2k+3},\ldots,x_\ell)$, $n'=n/p$ and $A'=U(n')^2$, then $S_1'$, $S_2'$ and $S_3'$ are $A'$-extremal sequences in $\Z_{n'}$ where $S_1'$, $S_2'$ and $S_3'$ denote the sequences in $\Z_n'$ which are obtained by dividing the terms of $S_1$, $S_2$ and $S_3$ by $p$. 
\end{theorem}

\begin{proof}
As $S$ is an $A$-extremal sequence in $\Z_n$, we have $\ell=3^{\Omega(n)}-1$. Suppose for any prime divisor $q$ of $n$ there are at least three terms of $S$ which are not divisible by $q$. Let $q$ be a prime divisor of $n$, $q^r\mid\mid n$ and let $A_q=U(q^r)$. Thus $S^{(q)}$ has at least three units where $S^{(q)}$ is as defined before Observation \ref{observation}. By Corollary \ref{cm} as $q$ is at least 7, it follows that $S^{(q)}$ is an $A_q^2$-weighted zero-sum sequence in $\Z_{q^r}$. 

As this is true for every prime divisor $q$ of $n$, by Observation \ref{observation2} we see that $S$ is an $A$-weighted zero-sum sequence. This contradicts our assumption about $S$. So we see that there is a prime divisor $p$ of $n$ such that at most two terms of $S$ are not divisible by $p$. Now by using Corollary \ref{g} we see that exactly two terms of $S$, say $x^*$ and $x^{**}$, are not divisible by $p$.

Let us assume that $x^*$ occurs before $x^{**}$ in $S$. Suppose $x^*\neq x_{k+1}$ or $x^{**}\neq x_{2(k+1)}$ where $k+1=(\ell+1)/3$. Hence, there is a subsequence $T$ of consecutive terms of $S$ of length at least $k+1$ such that $p$ divides every term of $T$. As every prime divisor of $n$ is at least 7 and as $T$ has length at least $(\ell+1)/3$, by Lemma \ref{g'} we see that $S$ has an $A$-weighted zero-sum subsequence of consecutive terms. This contradicts our assumption about $S$. Thus $x^*=x_{k+1}$ and $x^{**}=x_{2(k+1)}$.

Let $n'=n/p$ and for $1\leq i\leq 3$, let $S_i$ and $S_i'$ be defined as in the statement of the theorem. If $S_1'$ has an $A'$-weighted zero-sum subsequence of consecutive terms, then by Lemma \ref{c} we see that $S_1$ has an $A$-weighted zero-sum subsequence of consecutive terms. As $S_1$ is a subsequence of consecutive terms of $S$, it follows that $S$ has an $A$-weighted zero-sum subsequence of consecutive terms. This contradicts our assumption about $S$.

As $k+1=(\ell+1)/3$ and as $\ell+1=3^{\Omega(n)}$, it follows that $k+1=3^{\Omega(n')}$. Hence, as $S_1'$ has length $k$, it follows that $S_1'$ has to be an $A'$-extremal sequence in $\Z_{n'}$. Similar arguments show that $S_2'$ and $S_3'$ are also $A'$-extremal sequences in $\Z_{n'}$.
\end{proof}

\begin{lemma}\label{f}
Let $m$ be a divisor of $n$ and let $b\in U(m)$. Then there exists $a\in U(n)$ such that $a$ maps to $b$ under the natural map $\Z_n\to\Z_m$. If $b\in U(m)^2$, then there exists $a\in U(n)^2$ such that $a$ maps to $b$ under the natural map $\Z_n\to\Z_m$.
\end{lemma}

\begin{proof}
Let $\Z^+=\{i\in\Z : i>0\}$. As $b$ is coprime to $m$, by Dirichlet's theorem (see \cite{S}, for instance) there are infinitely many primes in the arithmetic progression $\{\,b+im : i\in\Z^+\,\}$. As only finitely many of these primes will divide $n$, we see that there exists $i\in\Z^+$ such that $a=b+im\in U(n)$. Also $a$ maps to $b$ under the natural map $\Z_n\to\Z_m$. The proof when $b\in U(m)^2$ now follows from this. 
\end{proof}

\begin{theorem}\label{8}
We continue with the notations which were used in the statement and proof of Theorem \ref{7}. Then the image of the sequence $(x^*,x^{**})$ under the natural map $\Z_n\to\Z_p$ is a $Q_p$-extremal sequence in $\Z_p$.
\end{theorem}

\begin{proof}
If $n=p$, then $S$ is a $Q_p$-extremal sequence in $\Z_p$. So we can assume that $n$ is not a prime. Suppose the image of the sequence $(x^*,x^{**})$ under the natural map $\Z_n\to\Z_p$ is a $Q_p$-weighted zero-sum sequence. Hence, by Lemma \ref{f} there exist $a,b\in U(n)^2$ such that $ax^*+bx^{**}\in\Z_n$ maps to $0\in\Z_p$ under the natural map $\Z_n\to\Z_p$. Thus $p$ divides $y=ax^*+bx^{**}$ in $\Z_n$. Let $T$ be the sequence which is obtained by removing the terms $x^*$ and $x^{**}$ from $S$ and by adding the term $y$ at the end of $S$. Now by using Theorem \ref{7} we see that $p$ divides all the terms of $T$.

Let $n'=n/p$ and let $A'=U(n')^2$. Consider the sequence $T'$ in $\Z_{n'}$ whose terms are obtained by dividing the terms of $T$ by $p$. Let $q$ be a prime divisor of $n'$. As $S_1'$ is an $A'$-extremal sequence in $\Z_{n'}$ and as every prime divisor of $n'$ is at least 7, by Corollary \ref{g} at least two terms of $S_1'$ are coprime to $q$. Similarly, at least two terms of $S_2'$ are coprime to $q$. Hence, $T'$ has at least 4 terms which are coprime to $q$. As this is true for any prime divisor $q$ of $n'$, by the same argument as in the first paragraph of the proof of  Theorem \ref{7} we see that $T'$ is an $A'$-weighted zero-sum sequence in $\Z_{n'}$. Now by Lemma \ref{c} we see that $T$ is an $A$-weighted zero-sum sequence in $\Z_n$.

As $T$ is a concatenation of the sequences $S_1,~S_2,~S_3$ and $(y)$, and as $a,b\in A$, it follows that $S$ is an $A$-weighted zero-sum sequence in $\Z_n$. This contradicts our assumption about $S$. Thus, the image of the sequence $(x^*,x^{**})$ under the natural map $\Z_n\to\Z_p$ cannot be a $Q_p$-weighted zero-sum sequence in $\Z_p$. As $C_{Q_p}(p)=3$ for $p>2$, it follows that the image of this sequence is a $Q_p$-extremal sequence in $\Z_p$.
\end{proof}

\begin{remark}
Theorems \ref{sq}, \ref{7} and \ref{8} together characterize the $U(n)^2$-extremal sequences in $\Z_n$ where any prime divisor of $n$ is at least $7$. 
\end{remark}

We conclude this section with the following results about these sequences. 

\begin{proposition}
Using the notation of Theorem \ref{7}, we have that p is the unique prime divisor of $n$ such that $p$ is coprime to exactly two terms of $S$, and any other prime divisor $q$ of $n$ is coprime to at least three terms of $S$.
\end{proposition}

\begin{proof}
By Theorem \ref{7} we know that $p$ is coprime to exactly two terms of $S$. Let $q$ be a prime divisor of $n$ other than $p$. By Corollary \ref{g} we see that $q$ has to be coprime to at least two terms of $S$. We claim that $q$ cannot be coprime to exactly two terms of $S$. Suppose our claim is false. As a result, by a similar argument as in the second paragraph of the proof of Theorem \ref{7} we see that those two terms have to be $x^*$ and $x^{**}$. So $q$ divides all the terms of $S_1$. As $q\neq p$, this implies that $q$ divides all the terms of $S_1'$. However, as $q$ is a divisor of $n'$ and as $S_1'$ is an $A'$-extremal sequence, by Corollary \ref{g} we see that $q$ has to be coprime to at least two terms of $S_1'$. Thus we get a contradiction. Hence, our claim must be true.  
\end{proof}

\begin{proposition}\label{obsu2}
Let $A=U(n)^2$ where $n=p_1p_2$ is a product of two not necessarily distinct primes which are at least 7. Then  a sequence $S$ in $\Z_n$ is an $A$-extremal sequence if and only if $S$ is of the form 
$(b_1q_1,\,b_2q_1,\,a_1,\,b_3q_1,\,b_4q_1,\,a_2,\,b_5q_1,\,b_6q_1)$ 
where $q_1,q_2$ is a permutation of $p_1,p_2$, the sequence $(b_i,\,b_{i+1})$ is a $Q_{q_{_2}}$- extremal sequence in $\Z_{q_{_2}}$ for $i=1,3,5$ and the image of the sequence $(a_1,a_2)$ under the natural map $\Z_n\to\Z_{q_{_1}}$ is a $Q_{q_{_1}}$- extremal sequence in $\Z_{q_{_1}}$.   
\end{proposition}

\begin{proof}
As $S$ is an $A$-extremal sequence in $\Z_n$, it has length $8$. Let $S$ be the sequence $(x_1,\ldots,x_8)$. By Theorem \ref{7} there is a prime $q_1\in\{p_1,p_2\}$ such that $x_3,x_6$ are the only terms of $S$ which are coprime to $q_1$. Let $n'=n/q_1$, $A'= U(n')^2$ and let $S_1'=(x_1/q_1,x_2/q_1)$, $S_2'=(x_4/q_1,x_5/q_1)$, $S_3'=(x_7/q_1,x_8/q_1)$.

By Theorem \ref{7} we get that $S_1'$, $S_2'$ and $S_3'$ are $A'$-extremal sequences in $\Z_{n'}$ and so they are $Q_{q_{_2}}$- extremal sequences in $\Z_{q_{_2}}$ where $q_2$ is the prime $n'$. By Theorem \ref{8} we get that the image of the sequence $(x_3,x_6)$ under the natural map $\Z_n\to\Z_{q_{_1}}$ is a $Q_{q_{_1}}$- extremal sequence in $\Z_{q_{_1}}$. Hence, we see that $S$ is of the form as mentioned in the statement of Proposition \ref{obsu2}. Now Theorem \ref{sq} shows that any such sequence will be an $A$-extremal sequence in $\Z_n$.  
\end{proof}

\section{The case when $A=U(n)^3$}\label{cubes} 

Let $n=p_1p_2\ldots p_s$ be squarefree. Let $n_1$ be the divisor of $n$ such that for any prime divisor $p$ of $n$, we have $p\mid n_1$ if and only if $p\equiv 1~(\textrm{mod}~3)$. Also let $n_2=n/n_1$. We will follow the {\bf notation} $n=n_1n_2$ for the rest of this section. 

The next result follows from Lemmas 1 and 6 of $\cite{SKS}$, by using the fact that if $p\equiv 2~(\textrm{mod}~ 3)$, then $U(p)^3=U(p)$. 

\begin{lemma}[\cite{SKS}]\label{cub5}
Let $A=U(p)^3$ where $p$ is a prime such that $p\neq 2,7,13$. Then a sequence $S$ in $\Z_p$ is an $A$-weighted zero-sum sequence if at least three elements of $S$ are in $U(p)$. 
\end{lemma}

\begin{remark}
We can check that the sequence $(1,1,1)$ in $\Z_p$ is not a $U(p)^3$-weighted zero-sum sequence when $p=2,7,13$. So the conclusion of Lemma \ref{cub5} is not true for these values of $p$ (as was mentioned in the Remark after Lemma 1 of \cite{SKS}). 
\end{remark}

In Corollary 8 of $\cite{SKS}$ it was shown that for $A=U(n)^3$ where $n$ is squarefree and is not divisible by 2, 7 or 13, we have $C_A(n)=2^{\Omega(n_2)}3^{\Omega(n_1)}$. We now give a method to construct $A$-extremal sequences in $\Z_n$.

\begin{theorem}\label{cub6}
Let $A=U(n)^3$ where $n$ is squarefree and is not divisible by 2,7 or 13. Let $p$ be a prime divisor of $n$, $n'=n/p$ and $A'=U(n')^3$. Then the following constructions give us $A$-extremal sequences in $\Z_n$ in each of the cases below.

If $p\equiv 1~(\textrm{mod}~3)$, let $S=(p\,u_1,\,\ldots,\,p\,u_k,\, x^*,\,p\,v_1,\,\ldots,\,p\,v_k,\,x^{**},\, p\,w_1,\,\ldots,\,p\,w_k)$ where $x^*,x^{**}\in\Z_n$ are such that the image of the sequence $(x^*,x^{**})$ under the natural map $\Z_n\to\Z_p$ does not have any $U(p)^3$-weighted zero-sum subsequence and $S_1'=(u_1,\,\ldots,\,u_k)$, $S_2'=(v_1,\,\ldots,\,v_k)$, $S_3'=(w_1,\,\ldots,\,w_k)$ are $A'$-extremal sequences in $\Z_{n'}$.

If $p\not\equiv 1~(\textrm{mod}~3)$, let $S=(p\,u_1,\,\ldots,\,p\,u_k,\, x^*,\,p\,v_1,\,\ldots,\,p\,v_k)$  where $x^*\in\Z_n$ is not divisible by $p$ and $S_1'=(u_1,\,\ldots,\,u_k)$, $S_2'=(v_1,\,\ldots,\,v_k)$ are $A'$-extremal sequences in $\Z_{n'}$. 
\end{theorem}

\begin{proof}
The proof of Theorem \ref{cub6} is very similar to the proofs of Theorem \ref{unit1'} and Theorem \ref{sq}. 
\end{proof}

Let us see two examples of how we can use Theorem \ref{cub6} to construct $U(n)^3$-extremal sequences.

As $U(5)^3=U(5)$, by Theorem \ref{nz} we have the $U(5)^3$-extremal sequences $(2),(2)$ and $(4)$ in $\Z_5$. Also by Corollary \ref{cubpx} we have that $(-1,2)$ is a $U(19)^3$-extremal sequence in $\Z_{19}$. Hence, by Theorem \ref{cub6} we get the $U(95)^3$- extremal sequence $(38,37,38,78,76)$ in $\Z_{95}$.

By Corollary \ref{cubpx} we have that $(3,6)$, $(5,7)$ are $U(19)^3$-extremal sequences in $\Z_{19}$. So from Theorem \ref{cub6} we get the $U(95)^3$- extremal sequence $(15,30,69,25,35)$ in $\Z_{95}$.

We can prove a result similar to Lemma \ref{unit4'} and Lemma \ref{g'} for $A=U(n)^3$, whose proof is also similar to the proofs of those lemmas.

\begin{lemma}\label{cub7}
Let $A=U(n)^3$ where $n$ is squarefree and is not divisible by 2, 7 or 13 and let $\ell=2^{\Omega(n_2)}3^{\Omega(n_1)}-1$. Let $S$ be a sequence in $\Z_n$ and let $q$ be a prime divisor of $n$. Suppose when $q$ divides $n_1$ there exists a subsequence $T$ of consecutive terms of $S$ of length at least $(\ell+1)/3$ such that each term of $T$ is divisible by $q$, and when $q$ divides $n_2$ there exists a subsequence $T$ of consecutive terms of $S$ of length at least $(\ell+1)/2$ such that each term of $T$ is divisible by $q$. Then $S$ has an $A$-weighted zero-sum subsequence of consecutive terms.
\end{lemma}

The proof of the next corollary is very similar to the proofs of Corollary \ref{unit4} and Corollary \ref{g}.  

\begin{corollary}\label{cub8}
Let $A=U(n)^3$ where $n$ is squarefree and is not divisible by 2, 7 or 13 and let $S$ be an $A$-extremal sequence in $\Z_n$. Then for any prime divisor $q$ of $n_1$ (resp. of $n_2$) we have that $q$ is coprime to at least two terms (resp. at least one term) of $S$. 
\end{corollary}

\begin{theorem}\label{cub9}
Let $A=U(n)^3$ where $n$ is squarefree and is not divisible by 2, 7 or 13 and let $S=(x_1,\ldots,x_\ell)$ be an $A$-extremal sequence in $\Z_n$. Then $S$ can be of one of the following two types.

Type 1: There is a prime divisor $p$ of $n_1$ such that $p$ divides all the terms of $S$ except the terms $x_{k+1}$ and $x_{2(k+1)}$ where $k+1=(\ell+1)/3$, and the image of the sequence $(x_{k+1},x_{2k+2})$ under the natural map $\Z_n\to\Z_p$ is a $U(p)^3$-extremal sequence in $\Z_p$. Moreover, if $S_1=(x_1,\ldots,x_k)$, $S_2=(x_{k+2},\ldots,x_{2k+1})$, $S_3=(x_{2k+3},\ldots,x_\ell)$, $n'=n/p$ and $A'=U(n')^3$, then $S_1'$, $S_2'$ and $S_3'$ are $A'$-extremal sequences in $\Z_{n'}$ where $S_1'$, $S_2'$ and $S_3'$ denote the sequences in $\Z_n'$ which are obtained by dividing the terms of $S_1$, $S_2$ and $S_3$ by $p$.

Type 2: There is a prime divisor $p$ of $n_2$ such that $p$ divides all the terms of $S$ except the term $x_{k+1}$ where $k+1=(\ell+1)/2$. Moreover, if $S_1=(x_1,\ldots,x_k)$, $S_2=(x_{k+2},\ldots,x_\ell)$, $n'=n/p$ and $A'=U(n')$, then $S_1'$ and $S_2'$ are $A'$-extremal sequences in $\Z_{n'}$ where $S_1'$ and $S_2'$ denote the sequences in $\Z_{n'}$ which are obtained by dividing the terms of $S_1$ and $S_2$ by $p$.
\end{theorem}

\begin{proof}
As $S$ is an $A$-extremal sequence in $\Z_n$, it follows that $\ell=2^{\Omega(n_2)}3^{\Omega(n_1)}-1$. We consider the following three cases.

{\it Case (a): For any prime divisor $q$ of $n_1$ at least three terms of $S$ are coprime to $q$, and for any prime divisor $q$ of $n_2$ at least two terms of $S$ are coprime to $q$}.

In this case as $q\neq 2,7,13$, if $q$ divides $n_1$, then by Lemma \ref{cub5} we have that $S^{(q)}$ is a $U(q)^3$-weighted zero-sum sequence in $\Z_{q}$. As $q$ is an odd prime, if $q$ divides $n_2$, then by Lemma \ref{gri} we see that $S^{(q)}$ is a $U(q)$-weighted zero-sum sequence in $\Z_{q}$.

So for any prime divisor $q$ of $n$ we get that $S^{(q)}$ is a $U(q)^3$-weighted zero-sum sequence in $\Z_{q}$ and hence by Observation \ref{observation2} we see that $S$ is an $A$-weighted zero-sum sequence. This contradicts our assumption about $S$.

{\it Case (b): There is a prime divisor $p$ of $n_1$ such that at most two terms of $S$ are coprime to $p$}.

In this case the proof of the result is very similar to the proofs of Theorems \ref{7} and \ref{8} where we use Lemma \ref{cub5}, Lemma \ref{cub7} and Corollary \ref{cub8} in place of Corollary \ref{cm}, Lemma \ref{g'} and Corollary \ref{g}. We also use a consequence of Lemma \ref{f}.

{\it Case (c): There is a prime divisor $p$ of $n_2$ such that at most one term of $S$ is coprime to $p$}.

In this case the proof of the result is very similar to the proof of Theorem \ref{unit3} where we use Lemma \ref{cub7} and Corollary \ref{cub8} in place of Lemma \ref{unit4'} and Corollary \ref{unit4}.
\end{proof}

\begin{remark}
Theorems \ref{cub6} and \ref{cub9} together characterize the $U(n)^3$-extremal sequences in $\Z_n$ where $n$ is squarefree and is not divisible by 2, 7 or 13. 
\end{remark}

\section{Concluding remarks}

In Theorem \ref{unit3} we have classified the $U(n)$-weighted extremal sequences in $\Z_n$ when $n$ is odd. The corresponding result for an even integer could be studied. In the characterization of $U(n)^3$-weighted extremal sequences in $\Z_n$, it will be interesting if we could remove the assumption that $n$ is squarefree.

In Lemma \ref{f} we have proved that the image of $U(n)$ under the natural map $\Z_n\to\Z_m$ is $U(m)$ using Dirichlet's theorem. In \cite{CM} an elementary proof of this fact is given just before Section 3.

In Definition \ref{eq} if a sequence $T$ is obtained by a permutation of the terms of a sequence $S$, then we do not say that $T$ is $A$-equivalent to $S$. This is because if $S$ has an $A$-weighted zero-sum subsequence of consecutive terms, then the same may not be true for a sequence $T$ which is obtained in this manner. This is why Definition \ref{eq} differs from the definition which is given at the end of Section 2 in \cite{AMP}.

\bigskip

\noindent {\bf Acknowledgement.}
We sincerely thank the editor for taking the time to  provide suggestions which have helped to improve the presentation of this paper. Santanu Mondal would like to acknowledge CSIR, Govt. of India for a research fellowship.


\begin{thebibliography}{1}\footnotesize
\bibitem{AC} S. D. Adhikari and Y. G. Chen, Davenport constant with weights and some related questions, II, {\em J. Combin. Theory Ser. A} {\bf 115}, no. 1, (2008), 178--184. 

\bibitem{ACF}
S. D. Adhikari, Y. G. Chen, J. B. Friedlander, S. V. Konyagin and F. Pappalardi, Contributions to zero-sum problems, {\em Discrete Math.} {\bf 306} (2006), 1--10.

\bibitem{AMP} S. D. Adhikari, I. Molla and S. Paul, Extremal sequences for some weighted zero-sum constants for cyclic groups, {\em CANT IV, Springer Proc. Math. Stat.} {\bf 347} (2021), 1--10.

\bibitem{AR}
S. D. Adhikari and P. Rath,
Davenport constant with weights and some related questions, {\em Integers} {\bf 6} (2006), A30. 

\bibitem{CM} M. N. Chintamani and B. K. Moriya, Generalizations of some zero sum theorems, {\em Proc. Indian Acad. Sci. Math. Sci.} {\bf 122}, no. 1, (2012), 15--21. 

\bibitem{sg} Simon Griffiths, The Erd\H{o}s-Ginzberg-Ziv theorem with units, {\em Discrete Math.} {\bf 308}, no. 23, (2008), 5473--5484. 

\bibitem{SKS} S. Mondal, K. Paul and S. Paul, On a different weighted zero-sum constant, {\em Discrete Math.}, to appear.  http://arxiv.org/abs/2110.02539

\bibitem{S} Atle Selberg, An elementary proof of Dirichlet's theorem about primes in an arithmetic progression, {\em Ann. of Math. (2)} {\bf 50} (1949), 297--304. 
\end{thebibliography}
\end{document}